\documentstyle[graphicx]{article}
\author{Vassily Olegovich Manturov}

\title{Minimal diagrams of classical and virtual links}

\newtheorem{thm}{Theorem}
\newtheorem{lm}{Lemma}
\newtheorem{re}{Remark}

\date{}

\begin{document}

\maketitle{To Lou Kauffman, on the occasion of his 60th birthday.}

\abstract{We prove that a virtual link diagrams satisfying two
conditions on the Khovanov homology is minimal, that is, there is
no virtual diagram representing the same link with smaller number
of crossings. This approach works for both classical and virtual
links}

\vspace{2cm}

For definitions of the Jones polynomial, Kauffman bracket, and the
Khovanov homology, we send the reader to the original papers
\cite{Jon,Kau,Kh}.

The theory of virtual links was first proposed by Kauffman in
1996, see \cite{KaV}.

We shall treat virtual links as a combinatorial generalisation of
classical links with a new crossing type, called {virtual},
allowed. The set of Reidemeister moves is enlarged by a new move,
called {\em the detour move}, which means the following. Having a
virtual diagram $L$ with an arc $AB$ consisting of some
consecutive virtual crossings (and no classical crossings), one
can remove this arc and draw it elsewhere as a curve connecting
the same endpoints, $A$ and $B$, so that all new crossing are set
to be virtual. Thus, a virtual link is an equivalence class of
virtual link diagrams modulo Reidemeister moves and detour moves,
for more details, see \cite{KaV}.

With any virtual link diagram $L$, one associates an {\em atom},
i.e. two-dimensional surface together with a $4$-graph embedded
into and a checkerboard colouring of $2$-cells.

The vertices of the atom come from the {\em classical crossings}
of the diagrams, whence the rules for attaching the black cells
come from over/ undercrossing structures. We decree ``black
angles'' to be those lying between overcrossing and undercrossing
such that while moving from an undercrossing to an overcrossing
clockwise, we sweep the black angle.

One recovers the initial diagram from the atom up to
virtualisations and detours. Denote the atom corresponding to the
link $L$ by $V(L)$, its Euler characteristic by $\chi(L)$, and its
genus by $g(L)$. Note that the atom need not be orientable. Thus,
e.g., the atom corresponding to the virtual trefoil lives on the
2-dimensional projective plane. Also, the genus of the atom does
not coincide with the underlying genus of thickened surface where
the virtual knot lives, see, e.g. \cite{Kup}. Indeed, the
underlying surface need not admit any checkerboard colouring. By a
virtualisation (see,e.g., \cite{Ma}) we mean a replacement of a
classical crossing by a (2-2)-tangle consisting of three
consecutive crossings: a virtual one, the classical one (with the
same writhe as the initial one) and one virtual crossing. The atom
does not feel the virtualisation of the initial diagram whence the
checkerboard surface does.

All necessary constructions concerning knots and atoms can be
found in \cite{Ma}.

By span of a 1-variable Laurent polynomial we mean the difference
between its leading degree and its lowest degree.

\begin{lm}[\cite{Ma}]. Given a virtual diagram $L$. Then $span\langle L\rangle \le
4n+2(\chi(L)-2)$.
\end{lm}

Here $n$ is the number of classical vertices of the diagram $L$.

Indeed, this quantity $4n+2(\chi(L)-2)$ appears while considering
the A-state and the B-state of the Kauffman state sum expansion.
We should just take into account that (by definition) the number
of white cells of the atom equals the number of curves in the
A-state, and the number of black cells of the atom equals the
number of curves in the B-state.

Let us say that a virtual diagram $L$ satisfies {\em the first
completeness condition} if the inequality in Lemma 1 becomes a
strict equality. In other words, we have the first completeness
condition when neither the leading nor the lowest coefficient of
the bracket polynomial vanishes.

A very important question is to classify all those links having a
diagram satisfying the first completeness condition (briefly,
1-complete).

One can slightly generalise the first completeness equality in
different ways. For instance, one can use the span of the
$\Xi$-polynomial (see \cite{MaXi}) in the case of virtual knots
and take span in the sense of leading and lowest degree of its
coefficients. This polynomial has the same state sum expansions as
the Kauffman bracket, the main difference being some new
``geometric coefficients'' at monomials. Thus, $span\Xi(L) \le
4n+2(\chi(L)-2)$. In this case, the span might be larger than that
of the Kauffman bracket. Also, using the Khovanov polynomail
$Kh(q,t)$, one can obtain the bracket polynomial after a variable
change $t=-1$. Thus, it readily follows from the above discussion
that $span_{q}(Kh(q,t))\le 2n+\chi(L)$.

The Khovanov homology is well defined over arbitrary field of
coefficients for the case of classical links. Note that the
Khovanov homology can be used for the $Z_{2}$ case for arbitrary
virtual links link diagrams and with arbitrary field coefficients
for oriented link diagrams (in the sense of atoms), \cite{MaKho}.

Thus, we get some conditions slightly weaker than the first
completeness conditions. We say that a virtual link diagram is
{\em 1-complete in a broad sense} if either $span
\Xi(L)=4n+2(\chi(L)-2)$ or $span_{q}(K(q,t))\le 2n+\chi(L)$ for
Khovanov homologies over some ring of coefficients.

This leads us to examples where the first completeness conditions
fails. Thus, for instance, there exist a link $L$ for which the
leading term of the Khovanov polynomial (with respect to $q$) is
some polynomial $P(t)$ such that $P(-1)=0$, whereas $P(-1)$
coincides with the leading (or lowest, which does not matter)
coefficient of the Kauffman polynomial.

It would be very interesting to classify all those (classical and
virtual) links not admitting 1-complete diagrams in the broad
sense.

This problem remains actual if we allow the virtualisation move:
this move does not change the Kauffman bracket, it does not change
the Khovanov homology either. Thus we get a weaker equivalence on
virtual knot diagrams, which is definitely worth studying,
especially, in view of minimality problems.

From the above discussion, we obtain the following

\begin{thm} If a virtual link diagram $L'$ is 1-complete (in a broad sense) then
we can decrease the number of its vertices only at the expense of
the genus. \end{thm}

In other words, if a diagram $L'$ (of a link $L$) has $n$
classical vertices and genus $g$ and $L$ has a diagram $L''$ with
$n'<n$ classical vertices then $g(L'')<g(L')$.

An immediate corollary is the Kauffman-Murasugi theorem on
alternating links: they live on the sphere, thus having maximal
possible Euler characteristic equal to $2$ (resp., minimal
possible genus). In this case, the 1-completeness yields the
minimality.

The question now is {\bf how to handle the genus}?

From now on, for the sake of simplicity, we deal only with
orientable link diagrams (in the sense of atoms). On one hand, the
case of non-orientable diagrams can be handled by means of
double-coverings \cite{MaKho}. On the other hand, all arguments
below work in the unorientable case, see Remark \ref{rema}.

Here we have some ``grading'' on the set of (orientable) virtual
link diagrams: the lowest level is represented by alternating
classical knots and quasialternating knots (obtained from former
ones by virtualisations and detours), the next levels are
regulated by the genus of the corresponding atom.

It is conceptually important that classical knots should not be
considered separately from virtual knots; all results work in both
categories.

We shall use the following result.

\begin{thm} Let $L$ be a virtual link diagram. Then the
(unnormalised) Khovanov complex for $L$ is quasi-isomorphic to the
complex of the form $\sum_{s\in
K_{1}(L)}A[r(s)]\{r(s)\}[w(L_{s})]\{2w(L_{s})\}$, where $A$ is the
two-term complex with terms $v_{\pm}$ of grading $(0,\pm 1)$.
\end{thm}

Herewith, we use the spanning tree model for the Khovanov
homology. Here $K_{1}$ is the set of states each of which having
precisely one circle. For any $s$ we take some diagram $L_{s}$
obtained from $L$ by smoothing {\em some} crossings in such a way
that: we can get to the state $s$ by smoothing the remaining
circles of $L_{s}$ and the diagram $L_{s}$ can be unknotted only
by using the first Reidemeister move, for more details, see
\cite{Weh}. Also, $r(s)$ is the number of $B$-type smoothings of
the diagram $s$ and $w(s)$ is minus the writhe number (in Wehrli's
setting).

This theorem was proved by Stephan Wehrli \cite{Weh} for the case
of classical links. The proof in the virtual case belongs to the
author of the paper. Namely, one should take the construction from
\cite{MaKho} and repeat Wehrli's proof. The only thing to mention
is that in the case of {\em unorientable} virtual knot diagrams,
one can use only ${\bf Z}_{2}$-coefficients. However, we shall
deal only with {\bf oriented diagrams} and a {\bf field} of
coefficients. The proof of this generalisation is literally the
same; one should just accurately use the Khovanov homology
construction proposed in \cite{MaKho} and check all steps of the
proof.

\begin{re} As shown in my work \cite{MaKho}, Khovanov homology for orientable links
is invariant with coefficients in a given field (without torsions
in homology), the main difficulty being the K\"unneth formula for
the ``tensor square''. Perhaps, one can make a sharper statement
about arbitrary coefficients, but from now on, we deal only with
coefficients in a field.
\end{re}

The Khovanov homology lies on several diagonals, i.e. lines
$t-2q=const$. The thickness of the Khovanov homology is the number
of diagonals between the two extreme ones, i.e.,
$((t-2q)_{max}-(t-2q)_{min})/2+1$. Notation: $T(L)$

Let us prove the following important result.

\begin{thm} $T(L)\le g(L)+2$.\label{th3}
\end{thm}

Indeed, the diagonals correspond to the values of $r(s)$, i.e.,
numbers of $B$-smoothings in $1$-states. In the case of
alternating links (genus zero) it is known that $r(s)$ is the same
for all states $s\in K_{1}$. In the general case, we have to
estimate the amplitude of $r(s)$ for different $s$ from $K_{1}$.

\begin{re}
For unorientable link diagrams (and Khovanov homologies with ${\bf
Z}_{2}$-coefficients, we have indeed the same formula allowing the
genus $g(L)$ to be half-integer. All arguments remain the same;
the $q$-grading of the Khovanov homologies do not have the same
parity any more, for detailed description see \cite{MaKho}.
\label{rema}
\end{re}

The remaining part of Theorem \ref{th3} follows from

\begin{lm}. For an orientable link of genus $g$, the maximal possible
value of $r(s)$ and minimal possible value (for all states from
$K_{1}$) is equal to $2g$. Namely, it can be equal to $k,
k-2,\dots, k-2g$.
\end{lm}

The proof of this lemma goes as follows. We have some number of
curves $x$ in the $A$-state and some number of circles $y$ in the
$B$-state. To get to $K_{1}$ from the $A$-state, one should switch
at least $x-1$ crossings; also, to get to $K_{1}$ from the
$B$-state, one should switch at least $y-1$ crossings. Thus, the
values of $r(s)$ are in between $x-1$ and $n+1-y$. Now, one should
just recall the definition of the atom genus.

So, having two diagonals in Wehrli's complex for the case of
(quasi)alternating links, this number increases by one together
with the genus of the atom.

Thus, we know, where the {\bf chains} of Wehrli's complex live.
This tells us where to look for {\bf homologies} of Wehrli's
complex which coincide with Khovanov's homologies (unnormalised).
From this we deduce Theorem 3.

Having this, we define a(n orientable virtual) link to be {\bf
$2$-complete} if the number of diagonals is as large as it should
be, i.e., $g+2$, where $g$ is the genus of the link. In other
words, one should just take care that any of the two extreme
diagonals in Wehrli's complex have at least one homology.

From what above, we obtain the following

\begin{thm} (THE MINIMALITY THEOREM). Suppose an orientable virtual link
$L$ is 1-complete and 2-complete. Then this diagram is minimal.
\end{thm}

Indeed, from 2-completeness we see that we can not reduce the
genus. Together with 1-completeness this implies that we cannot
decrease the number of crossings.

\begin{re}
The minimality theorem remains true if we understand the
1-completeness condition in the broad sense.
\end{re}

\begin{re}
Note that the minimality theorem for classical link $L$ says that
if two completeness conditions hold then for the link $L$ there is
neither classical diagrams nor virtual diagrams with strictly
smaller number of classical crossings.

In the general case, the question whether a {\em minimal classical
diagram} is minimal in virtual category, is still open.
\end{re}

A very important question is to study the interrelation between
the atom genus (that we have used) and the underlying genus and
their minimalities. They are not the same. {\bf They become the
same if we admit virtualisation}.

A fair question to ask is {\bf how large is the set described in
minimality theorem?}

In my opinion, this question belongs either to philosophy or to
measure theory. What I can say at least is that it is {\bf rather
large} for virtual knots and {\bf wider than just the class of
alternating links} in the classical case.

Definitely, both 1-completeness and 2-completeness are two very
important properties to study, and I think, the set of minimal
diagram detected by the Minimality Theorem can be enlarged by some
things like coverings (which can make) or almost-complete diagrams
with some estimates of the leading coefficient, and so on.

To support the Minimality Theorem, let us consider the knot
$13n_{3663}$ from Shumakovitch's paper \cite{Shu}. Suppose we do
not know how this knot (taken from some table) looks like but now
only the ${\bf Q}$-Khovanov homology, see below.

\begin{equation}
\begin{array}{||c||c|c|c|c|c|c|c|c|c|c|c|c|c|c|c||}
\hline \hline & -6 & -5 & -4 & -3 & -2 & -1 & 0 & 1 & 2 & 3 & 4 &
5 & 6 & 7 \cr

 \hline 13 & & & & & & & & & & & & & & 1 \cr

\hline 11& & & & & & & & & & & & & &  \cr

\hline 9& & & & & & & & & & & & 1 & 1 & \cr

\hline 7& & & & & & & & & & 1 & 1 &  &  & \cr

\hline 5& & & & & & & & & 1 & &  &  &  & \cr

\hline 3& & & & & & & & & 2 & 1 &  &  &  & \cr

\hline 1& & & & & & & 2 & 1 & &  &  &  &  & \cr

\hline -1& & & & & 1 & 1 & 1 & 1 & &  &  &  &  & \cr

\hline -3& & & & &  & 1 &  &  & &  &  &  &  & \cr

\hline -5& & & & 1 & 1 &  &  &  & &  &  &  &  & \cr

\hline -7& & 1 & &  &  &  &  &  & &  &  &  &  & \cr

\hline -9& &  & &  &  &  &  &  & &  &  &  &  & \cr

\hline -11 & 1 &  & &  &  &  &  &  & &  &  &  &  & \cr

\hline \hline

\end{array}
\end{equation}

This information is sufficient to prove that this diagram is
minimal.

Indeed, it has 13 crossings and 4 diagonals. Thus, its genus
cannot be less than $2$, and the span of the Kauffman polynomial
can not exceed $52-8=44$, that is, the span of the Khovanov
homology (with respect to  $q$) cannot exceed 24. But it equals
$24=2\cdot 13+\chi=26-2$: it occupies places between $-11$ and
$+13$. So, this knot diagram is 1-complete and 2-complete. Thus,
the diagram is minimal by the minimality theorem.

\begin{figure}
\centering\includegraphics[width=200pt]{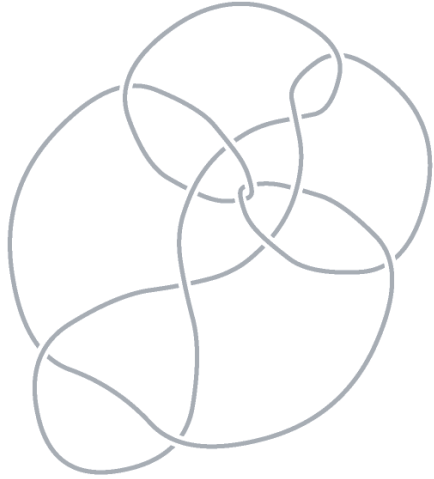} \caption{The
knot $13n_{3663}$}
\end{figure}

\end{document}